\documentclass[11pt,a4paper,twoside]{article}
\usepackage{amssymb}
\usepackage{amsmath}
\usepackage{fancyhdr}
\usepackage{amsthm}
\usepackage[cp1250]{inputenc}

\title{About the embedding of Moufang loops in alternative algebras}
\author{Nicolae Sandu}
\date{}
\begin{document}
\maketitle

\begin{abstract}

It is proved that any  free Moufang loop can be embedded in a loop
of invertible elements of some alternative algebra.
\smallskip\\

\noindent {Mathematics Subject Classification (1991):}  17D05,
20N05.
\smallskip\\
\noindent {Key words:} Moufang loop, alternative algebra.

\end{abstract}

\bigskip

Nowadays, the theory of loop algebras, in comparison with group
algebras, is not perfect at all. It is so because many loop laws
can't be transferred to loop algebras, while in group algebras the
group associativity  (and commuta\-tivity) is transferred to the
group algebras. So, if $L$ is  a  Moufang loop, then its loop
algebra $FL$ is not always alternative, i.e. the Moufang laws are
not always true in $FL$.

Till now the theory of loop algebras has generally developed when
the Moufang loops were  examined. For example, in [1 -  3] the
necessary and sufficient  conditions are pointed  for  the  loop
algebra $FL$ to be alternative, when $L$ is a Moufang  loop, and
in [1, 4 - 7] the construction of such loops $L$ is examined. We
note that these  conditions  are  quite  strong.  Loop  $L$ should
be ''almost'' associative.  Nowadays  the  theory  of such loops
is developing in such a way. It  is  assumed  beforehand that the
loop algebra $FL$ is alternative and there are examined either the
algebra $FL$, or Moufang loops $L$ setting different conditions on
algebra $FL$ [8 - 12]. These themes are stated in survey [13] and
[14] in details (where, in particular, the difficulties of this
theory are pointed).

It is well known that for an alternative algebra $A$ with  unit
the set $U(A)$ of all invertible elements of $A$ forms a Moufang
loop with respect to multiplication. This work offers another way
of examining Moufang  loops. Let $L$ be a free Moufang loop. It is
shown that if we factor the loop algebra $FL$ on some ideal $I$,
then $FL/I$ will be an alternative algebra and loop $L$ will be
embedded in the loop of invertible elements of algebra $FL/I$.
This is a positive answer to the question raised in [15]: is it
true that any Moufang loop can be imbedded into a homomorphic
image of a loop of type $U(A)$ for a suitable unital alternative
algebra $A$?  The equivalent version of this question is: whether
the variety generated by the loops of type $U(A)$ is a proper
subvariety of the variety of all Moufang loops?

The findings of this paper also give a partial positive answer to
a more general question (see, for example, [14]): is it true that
any Moufang loop can be imbedded into a loop of type $U(A)$ for a
suitable unital alternative algebra $A$? A positive answer to this
question was announced in [16]. Here, in fact, the answer to this
question is negative: in [15] it is constructed a Moufang loops
which are not imbedded into a loop of invertible elements of any
alternative algebra.

Let us now remember some notions and results from  the  loop
theory, which can be found in [13]. Loop $(L,\cdot) \equiv L$ is
called \textit{$IP$-loop} if the laws $^{-1}x\cdot xy = yx\cdot
x^{-1} = y$ are  true  in  it,  where $^{-1}xx = xx^{-1} = 1$. In
$IP$- loops $^{-1}x = x^{-1}$  and $(xy)^{-1} = y^{-1}x^{-1}$. The
loop is \textit{Moufang} if it satisfies the law
$$x(y \cdot zy) = (xy \cdot z)y. \eqno{(1)}$$
Every Moufang loop is a $IP$-loop. The subloop $H$ of loop $L$ is
called \textit{normal} in $L$, if
$$xH = Hx,\quad x\cdot yH = xy\cdot H, \quad H\cdot xy = Hx\cdot y \eqno{(2)}$$
for every $x, y \in L$.

Let $F$ be a field and $L$ be a loop. Let us examine the
\textit{loop algebra} $FL$. This is a free $F$-module with the
basis $\{q \vert q \in L\}$ and the product of the elements of
this basis is determined as their product in loop $L$. Let $H$ be
a normal subloop of loop $L$. We denote the ideal of algebra $FL$,
generated by the elements $1 - h$ ($h \in H$) by $\omega H$. If $H
= L$, then $\omega L$ is called the \textit{augmentation ideal} of
algebra $FL$. Let us determine the homomorphism of $F$-algebras
$\varphi$: $FL \to F(L/H)$ by the rule $\varphi(\sum \lambda_qq) =
\sum \lambda_qHq$. Takes place
\smallskip\\

\textbf{Lemma 1.} \textit{Let $H, H_1, H_2$ be normal subloops of
loop $L$. Then}

\textit{1) $\text{Ker}\varphi = \omega H$;}

\textit{2) $1 - h \in \omega H$ if and only if $h \in H$;}

\textit{3) if the elements $h_i$ generate the subloop $H$, then
the elements $1 - h_i$ generate the ideal $\omega H$; if $H_1 \neq
H_2$, then
 $\omega H_1 \neq \omega H_2$;
if $H_1 \subset H_2$, then $\omega H_1 \subset \omega H_2$; if $H
= \{H_1,H_2\}$, then $\omega H = \omega H_1 + \omega H_2$;}

\textit{4) $\omega L = \{\sum_{q \in L}\lambda_qq \vert \sum_{q
\in L}\lambda_q
 = 0\}$;}

\textit{5) $FL/\omega H \cong F(L/H),\quad \omega L/\omega H \cong
\omega (L/H)$;}

\textit{6) the augmentation ideal is generated as $F$-module by
the elements of the form $1 - q$ ($q \in L$).}
\smallskip\\

\textbf{Proof.}  As the mapping $\varphi$ is $F$-linear, then by
(2) for $h \in H, q \in L$ we have $\varphi((1 - h)q) = Hq - H(hq)
= Hq - Hq = 0$, i.e. $\omega H \subseteq \text{Ker}\varphi$. Let
now $K = \{k_j \vert j \in J\}$ be a complete system of
representatives of cosets of loop $L$ modulo the normal subloop
$H$ and let $\varphi(r) = 0$. We present $r$ as $r = u_1k_1 +
\ldots + r_tk_t$, where $u_i = \sum_{h \in H}\lambda_h^{(i)}h, k_i
\in K$. Then $0 = \varphi(r) = \varphi(u_1)\varphi(k_1) + \ldots +
\varphi(u_t)\varphi(k_t) = ( \sum_{h \in
H}\lambda_h^{(1)})\varphi(k_1) + \ldots + (\sum_{h \in
H}\lambda_h^{(t)})\varphi(k_t)$. As $\varphi(k_1), \ldots,
 \varphi(k_t)$ are pairwise distinct, then for all $i$ $\sum_{h \in
H}\lambda_h^{(i)} = 0$. Hence $-u_i = \sum_{h \in
H}\lambda_h^{(i)}(1 - h) - \sum_{h \in H}\lambda_h^{(i)} = \sum_{h
\in H}\lambda_h^{(i)}(1 - h)$ is an element from $\omega H$.
Consequently, $\text{Ker}\varphi \subseteq \omega H$, and then
$\text{ker}\varphi = \omega H$.

2). If $q \notin H$, then $Hq \neq H$. Then $\varphi(1 - q) = H -
Hq \neq 0$, i.e. by  1) $1 - q \notin \text{Ker}\varphi = \omega
H$.

3). Let elements $\{h_i\}$ generate subloop $H$ and $I$ be an
ideal, generated by the elements $\{1 - h_i\}$. Obviously $I
\subseteq \omega H$. Conversely, let $g \in H$ and $g = g_1g_2$,
where $g_1, g_2$ are words from $h_i$. We suppose that $1 - g_1, 1
- g_2 \in I$. Then $1 - g = (1-g_1)g_2 + 1 - g_2 \in I$,  i.e.
$\omega H \subseteq I$. Hence $I = \omega H$. Let $H_1 \neq H_2$
 (respect. $H_1 \subset H_2$) and $g \in H_1, g \notin H_2$. Then
 by    1) $1 - g \in \omega
H_1$, but $1 - g \notin \omega H_2$. Hence $\omega H_1 \neq \omega
H_2$ (respect. $\omega H_1 \subset \omega H_2)$. If $H = \{H_1,
H_2\}$, then by the first statement of  3) $\omega H = \omega H_1
+ \omega H_2$.

4). We denote $R = \{\sum_{q \in L}\lambda_qq \vert \sum_{q \in
 L}\lambda_q = 0\}$. Obviously, $\omega L \subseteq R$. Conversely,
 if $r \in R$ and $r = \sum_{q \in L}\lambda_qq$, then $-r =
- \sum_{q \in L}\lambda_qq = (\sum_{q \in L}\lambda_q)1 - \sum_{q
\in L}\lambda_qq = \sum_{q \in L}\lambda_q(1-q) \in \omega L$,
i.e. $R \subseteq \omega Q$. Hence $\omega L = R$.

5). Mapping $\varphi:FL \to F(L/H)$ is the homomorphism of loop
algebras and as by    1) $\text{Ker}\varphi = \omega H$, then
$FL/\omega H \cong F(L/H)$. Now from    4) it follows that $\omega
L/\omega H \cong \omega (L/H)$.

6). As $(1 - q)q^{\prime}= (1 - qq^{\prime}) - (1 - q^{\prime})$,
then the augmentation ideal $\omega L$ is generated by the
elements of form $1 - q$, where $q \in L$. This completes the
proof of Lemma 1.
\smallskip\\

Let $(B, +, \cdot)$ be an arbitrary algebra over a certain field
$F$. The mapping $\varphi$ of set $B$ is called the homomorphism
of algebra $(B, +, \cdot)$ if $\varphi(\lambda a) =
\lambda\varphi(a), \varphi(a + b) = \varphi a + \varphi b,
\varphi(a\cdot b) = \varphi a \cdot \varphi b$ for any $\varphi
\in F, a, b \in B$. If $(L,\cdot)$ is an arbitrary loop, $FL$ is
its loop algebra and $\varphi$ a certain homomorphism of algebra
$(FL, +, \cdot)$, then it follows from the last equality that  the
contraction $\varphi$ on $L$ will be the homomorphism of loop
$(L,\cdot)$. We call it \textit{$A$-homomorphism}.
$A$-homomorphism image $\varphi$ of loop $(L,\cdot)$ is not always
a loop, but only a groupoid with division. However it takes place.
\smallskip\\

\textbf{Lemma 2.} \textit{Let $(L,\cdot)$ be an $IP$-loop and let
$\varphi$ be a homomorphism of algebra $(FL, +, \cdot)$. Then
$A$-homomorphism image $\varphi$ of loop $(L,\cdot)$ will be a
loop.}
\smallskip\\

\textbf{Proof.} We denote the $A$-homomorphism image $\varphi$ of
loop $(L,\cdot)$ by $(\overline L,\star)$. It follows from the
$IP$-loop identity $x^{-1}\cdot xy = y$ that $\varphi(x^{-1}) =
(\varphi x)^{-1}$ and $(\varphi x^{-1})\star (\varphi
x\star\varphi y) = \varphi y, (\varphi x)^{-1}\star(\varphi
x\star\varphi y) = \varphi y, \overline x^{-1}\star(\overline
x\star \overline y) = \overline y$.  Let $\overline a, \overline b
\in \overline L$. It is obvious that the equation $\overline
a\star x = \overline b$ is always solvable and as $\overline
a^{-1}\star(\overline a\star x) = \overline a^{-1}\star \overline
b, x = \overline a^{-1}\star \overline b$, then it is uniquely
solvable. It can be shown by analogy that the equation $y\star
\overline a = \overline b$ is also uniquely solvable. Therefore,
$(\overline L,\star)$  is  a loop, as required.
\smallskip\\

Now, before we pass to the presentation of the basic  results,  we
give the construc\-tion of \textit{free $IP$-loop} with the set of
free generators $X  = \{x_1, x_2,\ldots\}$, using ideas from [13].
To the set $X$ we add the disjoint set
$\{x_1^{-1},x_2^{-1},\ldots\}$. Let us examine all groupoid words
$L(X)$ from set $\{x_1,x_1^{-1},x_2, x_2^{-1},\ldots\}$ relative
to multiplication $(\cdot)$ and let $e$ denote the empty word. For
the words from $L(X)$ we define the \textit{inverse words}: 1) for
$x_i \in X$ the inverse will be $x_i^{-1}$, and for $x_i^{-1}$ the
inverse will be $x_i$, i.e. ${(x_i^{-1})}^{-1} = x_i$;  2) if
$u\cdot v \in L(X)$, then $(u\cdot v)^{-1} = v^{-1}\cdot u^{-1}$.
Further, we define two words $u, v$ in $F(X)$ to be
\textit{Moufang-equivalent}, $u \approx v$ if one can be obtained
from other by a sequence of substitutions, each of which replaces
a subword $(rs\cdot r)t$ by $r(s\cdot rt)$ and inverse, where $r,
s, t$ are any words in $F(X)$. By a \textit{contraction} $\mu$ of
a word in $F(X)$ we mean the substitution a subword of the form
$u^{-1}(vw)$, $(wv)u^{-1}$, where $u \approx v$, by $w$. The
action $\nu$, inverse to contraction $\mu$ we call the
\textit{expansion}.

We define the $(\mu,\nu)$-equivalence $w \cong w^{\prime}$ for
words $w$, $w^{\prime}$ in $F(X)$ if one can be obtained from the
other by a finite sequence of substitutions each if which is
either a contradictions $\mu$ or expansion $\nu$ or a single use
of the Moufang law (1). The relation $\cong$ will be, obviously, a
relation of equivalence on $L(X)$. It will be even congruence, as
if a word $(u_1u_2\ldots u_n)_{\alpha}$ is given when $\alpha$ is
some parentheses distributions, obtained from words $u_1, u_2,
\ldots , u_n$, then the replacement of the word $u_i$, $i = 1, 2,
\ldots ,n$, with words or equivalence can be realized applying to
the given word a finite number of transformations of the above
described form.

With multiplication $\{u\}\cdot\{v\} = \{uv\}$ and inverse
$\{u\}^{-1} = \{u^{-1}\}$ of congruence classes we obtain a loop
with unity $\{e\}$, as the quotient loop $L(X)/\cong$ satisfies
the laws $x^{-1}\cdot xy = y$, $yx\cdot x^{-1} = y$. Moreover,
$L(X)/\cong$ will be a free Moufang loop on $\{x_i\}$, $i = 1, 2,
3, \ldots$ as set of free generators $X$. We identify $\{x_i\}$
with $x_i$ and we denote $L(X)/\cong$ by $L_X(\frak M)$.

Analogically of $F(X)$, we introduce the Moufang-equivalence,
transformations $\mu$, $\nu$ and $(\mu,\nu)$-equivalence for words
in $L_X(\frak M)$. We define a word in $L(\frak M)$ to be a
\textit{reduced word}  if no reductions of type $\mu$ of it are
possible.  If $w \in L(\frak M)$, then the number $l(w)$ of the
variables in $X$, contained in $w$, will be called the
\textit{length} of word $w$. Now let us show that if $w \to w_1, w
\to w_2$ are any reductions of type $\mu$  of a word $w$, then
there is a word $w_3$ obtained from each of $w_1, w_2$ by a
sequence of reductions of type look $\mu$. As a matter of fact, we
use induction on the length of $w$. If $l(w) = 1$, $w$ is already
a reduced word. If $l(w) = n$ and $w = u\cdot v$ where $u, v$ are
the subwords of $w$, then $l(u) < n, l(v) < n$. If both reductions
$w \to w_1, w \to w_2$ take place in the same subword, say $u$,
induction on length applied to $u$ yields the result. If the two
reductions take place in separate subwords, applying both gives
the $w_3$ needed. This leaves the case where at least one of the
reductions $w \to w_1, w \to w_2$ involves both subwords $u, v$ of
$w$. Then $w$ has, for example, the form $w = u^{-1}(uv)$.
Therefore $w = v$ and thus $l(w) < n$, then by inductive
hypothesis the statement is true.

Using this statement, one may prove by induction on length that
any word $w$ has a unique reduced words regarding the reductions
$\mu$. and all such reduced words belong to unique class of
Moufang-equivalence.Then, an induction on the number of reductions
and expansions connecting a pair of congruent words shows that
congruent words have the same reduced words.

Any word in $L_X(\frak M)$ has a reduced words. A \textit{normal
form} of a word $u$ in $L_X(\frak M)$ is a reduced word of the
least length. Clear, every word in $L_X(\frak M)$ has a normal
form. Let $u(x_1,x_2,\ldots , x_k)$, $u(y_1,y_2,\ldots ,y_n)$,
where $x_i, y_j \in X \cup X^{-1}$, be two words of normal form of
$u$ of length $l(u)$. $L_X(\frak M)$ is a free loop. Let, for
example, $y_1 \notin \{x_1,x_2,\ldots,x_k\}$, then
$u(x_1,x_2,\ldots,x_k) = u(1,y_2,\ldots,y_n)$. The length of
$u(1,y_2,\ldots,y_n)$ is strict least that $l(u)$. But this
contradicts the minimum condition for $l(u)$. Consequently, all
words of normal form of the same word in $L_X(\frak M)$ have the
same free generators in their structure. This completes the proof
of following statement.
\smallskip\\

\textbf{Lemma 3.} \textit{Any word in $L_X(\frak M)$ has a reduced
words that belongs to the unique class of Moufang-equivalence, two
words are $(\mu,\nu)$-equivalent if and only if they have the same
reduced words and all words of normal form of the same word in
$L_X(\frak M)$ have the same free generators in their structure.}
\smallskip\\

Now we consider a loop algebra $FM$ of free Moufang loop
$(M,\cdot) \equiv M$ over an arbitrary field $F$. Let $\overline M
= \{\overline u = 1 - u \vert u \in M\}$ and we define the
\textit{circle composition} $\overline u \circ \overline v =
\overline u + \overline v - \overline u \cdot \overline v$. Then
$(M, \circ )$ is a loop, denoted sometimes as $\overline M$. The
identity $\overline 1$ of $\overline M$ is the zero of $FM$,
$\overline 1 = 1 - 1$, and the inverse of $\overline u$ is
$\overline u^{-1} = 1 - u^{-1}$ as $\overline u \circ \overline 1
= 1 - u + 0 - (1 - u)0 = 1 - u = \overline u$, $\overline 1 \circ
\overline u = \overline u$, $\overline u \circ \overline u^{-1} =
\overline u + \overline u^{-1} - \overline u \overline u^{-1} = 1
- u + 1 - u^{-1} - (1 - u)(1 - u^{-1}) = 0$, $\overline u^{-1}
\circ \overline u = 0$. Let $\overline u, \overline v \in
\overline M$. Then $\overline u \circ \overline v = \overline u +
\overline v - \overline u\overline v = 1 - u + 1 - v - (1 - u)(1 -
v) = 1 -  uv = 1 - \overline{uv}$. Hence $\overline M$ is closed
under composition $(\circ)$ and
$$\overline u \circ \overline v = 1 - uy. \eqno{(3)}$$
Further, by (3) $\overline u^{-1}\circ(\overline u \circ \overline
v) = 1 - u^{-1}(uv) = 1 - v = \overline v$ and $(\overline v \circ
\overline u)\circ \overline u^{-1} = \overline v$. From here it
follows that $(\overline M,\circ)$ is a loop. We call it a
\textit{circle loop} corresponding to loop $(M,\cdot)$.

We define the one-to-one mapping $\overline \varphi: M \rightarrow
\overline M$ by $\overline \varphi(a) = \overline a$. For $a, b
\in M$ by (3) we have $\overline \varphi(ab) = 1 - ab = \overline
a \circ \overline b = \varphi(a) \circ \varphi(b)$. Hence
$\overline \varphi$ is an isomorphism of loop $M$ upon loop
$\overline M$. Then from Lemma 2 it follows that $\overline
\varphi$ induces the isomorphism $\varphi$ of loop algebra $FM$
upon loop algebra $F\overline M$ by rule $\varphi(\Sigma_{u\in
M}\alpha_uu) = \Sigma_{u\in M}\alpha_u(\overline \varphi(u)) =
\Sigma_{u\in M}\alpha_u\overline u$.

Clear, if the loop $M$ is generated by free generators $x_1, x_2,
\ldots$, then the loop $\overline M$ is generated by free
generators $\overline x_1, \overline x_2, \ldots$, the isomorphism
$\varphi: FM \rightarrow F\overline M$ is defined by mappings $x_i
\rightarrow \overline x_i$ and a word $u$ in $M$ has a normal form
if and only if the corresponding word $\overline u$ also has a
normal form. This completes the proof of following lemma.
\smallskip\\

\textbf{Lemma 4.} \textit{Let $FM$  be a loop algebra of  free
Moufang loop $(M,\cdot)$ with  free generators $x_1, x_2, \ldots$
and let $\overline M = \{\overline u = 1 - u \vert u \in M\}$ be
the corresponding loop under circle composition $\overline u \circ
\overline v = \overline u + \overline v - \overline u\overline v$.
Then the mappings $x_i \rightarrow \overline x_i$ define an
isomorphism $\varphi$ of loop algebra $FM$ upon loop algebra
$F\overline M$ by rule $\varphi(\Sigma\alpha_uu) =
\Sigma\alpha_u(\overline \varphi(u)) = \Sigma\alpha_u\overline u$,
$\alpha_u \in F$, $u \in M$, and a word in loop $(M,\cdot)$ has a
normal form if and only if the word $\varphi u$ has a normal form
in loop $(\overline M,\circ)$.}
\smallskip\\

Further, according to Lemma 4 for algebra $FM$ we will consider
only monomials of normal form. Let $u \in FM$ and let $\varphi$ be
the isomorphism defined in Lemma 4. We denote $\varphi(u) =
\overline u$. If $u = \Sigma\alpha_iu_i$, $\alpha_i \in F$, $u_i
\in M$, is a polynomial in $FM$ then we denote $\frak c(u) =
\Sigma\alpha_i$. Clear, that $\frak c(u) = \frak c(\overline u)$,
where $\overline u = \Sigma\alpha_i\overline u_i$.

If the free Moufang loop $M$ is non-associative, then from the
definition of loop algebra there follow the equalities
$$(a,b,c)+(b,a,c) = 0,\quad (a,b,c)+(a,c,b) = 0\quad \forall a, b, c \in L,
\eqno{(4)}$$ where the notation $(a,b,c) = ab\cdot c - a\cdot bc$
means that the associator in algebra does not always hold  in
algebra $FM$. Let $I(M)$ denote the ideal of algebra $FM$,
generated by all the elements of the left part of equalities (4).
It follows from the definition of  loop  algebra and
di-associativity of Moufang loops  that $FM/I(M)$ will be an
alternative algebra. We remind that algebra $A$ is called
\textit{alternative} if the identities $(x,x,y) = (y,x,x) = 0$
hold  in it. Hence we proved.
\smallskip\\

\textbf{Lemma 5.} \textit{Let $FM$ and $F\overline M$ be the loop
algebras of free Moufang loop $(M,\cdot)$ and its corresponding
circle loop $(\overline M,\circ)$ and let $I(M,\cdot)$,
$I(\overline M,\circ)$ be the ideals of $FM$ and $F\overline M$
respectively, defined above. Then $I(M) = I(\overline M)$ and for
any $\overline u \in I(\overline M)$ $\frak c(\overline u) = 0$.}
\smallskip\\

\textbf{Proof.} We denote $v_1 = v_1(u_{11},u_{12},u_{13}) =
(u_{11},u_{12},u_{13}) + (u_{12},u_{11},u_{13})$, $v_2 =
v_2(u_{21},u_{22},u_{23}) = (u_{21},u_{22},u_{23}) +
(u_{21},u_{23},u_{22})$, where $u_{ij} \in M$, $i = 1, 2$, $j = 1,
2, 3$. Then as $F$-module the ideal $I(M)$ is generated by the
elements of form
$$w(d_1,\ldots,d_k,v_i,d_{k+1},\ldots,d_m),$$
where $i = 1,2$ and $d_1,\ldots,d_m$ are monomials from $FM$.

Let $w = w(d_1,\ldots,d_k,v_1,d_{k+1},\ldots,d_m)$. Then by (3)
$$w = w(d_1,\ldots,d_k,(u_{11},u_{12},u_{13}) +
(u_{12},u_{11},u_{13}),d_{k+1},\ldots,d_m) =$$
$$w(d_1,\ldots,d_k,u_{11}u_{12}\cdot u_{13},d_{k+1},\ldots,d_m)
-$$
$$w(d_1,\ldots,d_k,u_{11}\cdot u_{12}u_{13},d_{k+1},\ldots,d_m)+$$
$$w(d_1,\ldots,d_k,u_{12}u_{11}\cdot u_{13},d_{k+1},\ldots,d_m)-$$
$$w(d_1,\ldots,d_k,u_{12}\cdot u_{11}u_{13},d_{k+1},\ldots,d_m)=$$
$$-(1-w(d_1,\ldots,d_k,u_{11}u_{12}\cdot
u_{13},d_{k+1},\ldots,d_m))+$$
$$(1-w(d_1,\ldots,d_k,u_{11}\cdot u_{12}u_{13},d_{k+1},\ldots,d_m))-$$
$$(1-w(d_1,\ldots,d_k,u_{12}u_{11}\cdot u_{13},d_{k+1},\ldots,d_m))+$$
$$(1-w(d_1,\ldots,d_k,u_{12}\cdot u_{11}u_{13},d_{k+1},\ldots,d_m)) =
$$
$$\overline w(\overline d_1,\ldots,\overline d_k,(\overline u_{11}\circ \overline
u_{12})\circ  \overline u_{13},\overline d_{k+1},\ldots,\overline
d_m) -$$
$$\overline w(\overline d_1,\ldots,\overline d_k,\overline u_{11}\circ (\overline
u_{12}\circ \overline u_{13}),\overline d_{k+1},\ldots,\overline
d_m)+$$
$$\overline w(\overline d_1,\ldots,\overline d_k,(\overline u_{12}\circ \overline
u_{11})\circ  \overline u_{13},\overline d_{k+1},\ldots,\overline
d_m)-$$
$$\overline w(\overline d_1,\ldots,\overline d_k,\overline u_{12}\circ
(\overline u_{11}\circ \overline u_{13}),\overline
d_{k+1},\ldots,\overline d_m) = $$
$$\overline w(\overline d_1,\ldots,\overline d_k,\overline v_2,\overline d_{k+1},\ldots,\overline
d_m).$$ Similarly, $w(d_1,\ldots,d_k,v_2,d_{k+1},\ldots,d_m) =
\overline w(\overline d_1,\ldots,\overline d_k,\overline
v_2,\overline d_{k+1},\ldots,\overline d_m)$. Hence $I(M)
\subseteq I(\overline M)$.

Conversely, we consider a polynomial in $f\overline M$ of form
$\overline w(\overline d_1,\ldots,\overline d_k,\overline
v_i,\break \overline d_{k+1},\ldots,\overline d_m)$. It is clear
that $\overline w \in I(\overline M)$ and any element $\overline z
\in I(\overline M)$ will be represented as sum of finite number of
polynomials of a such form. We have $\frak c(\overline v_i) = 0$,
then $\frak c(\overline w) = 0$ and, consequently, $\frak
c(\overline z) = 0$. Now, let for example $\overline v_i =
\overline v_1$. By (3) we get $\overline v_1 = (\overline
u_{11}\circ \overline u_{12})\circ \overline u_{13} = \overline
u_{11}\circ (\overline u_{12}\circ \overline u_{13}) = 1 -
u_{11}u_{12}\cdot u_{13} - (1 - u_{11}\cdot u_{12}u_{13}) =
-u_{11}u_{12}\cdot u_{13} + u_{11}\cdot u_{12}u_{13} =
-(u_{11},u_{12},u_{13}) = -v_1$. Further, by relation $\overline
x\circ \overline y = 1 -xy$ in a expression $\overline w$ we pass
from operation $(\circ)$ to operation $(\cdot)$. Then $\overline
w$ can be written as the sum of a finite number of monomials, each
of them containing the associators $v_i$ in its structure. Then
$\overline w \in I(M)$, and hence $\overline z \in I(M)$,
$I(\overline M) \subseteq I(M)$. Consequently, $I(\overline M) =
I(M)$. This completes the proof of Lemma 5.
\smallskip\\

\textbf{Theorem 1.} \textit{Let $(M,\cdot)$ be a free Moufang
loop, let $F$ be an arbitrary field and let $\varphi: FM
\rightarrow FM/I(M)$ be the natural homomorphism of algebra $FM$
upon the alternative algebra $FM/I(M)$. Then the image
$\varphi(M,\cdot) = (\overline M,\star)$ of loop $M,\cdot)$ will
be the isomorphism of these loops.}
\smallskip\\

\textbf{Proof.} Any Moufang loop is an $IP$-loop, then by Lemma 2
the image of loop $(M,\cdot)$ under the $A$-homomorphism $\varphi:
FM \rightarrow FM/I(M)$ will be a loop $(\overline M,\star)$. Let
$H$ be a normal subloop of loop $(M,\cdot)$, that corresponds to
$\varphi$. Then $1 - H \subseteq I(M)$. We suppose that $H \neq
\{1\}$ and let $1 \neq u(x_1,\ldots ,x_k) \in H$ be a word in the
free generators $x_1, \ldots , x_k$ of normal form. Then the
length $l(u) > 0$. By (3) we write $1 - u(x_1, \ldots , x_k)$ in
generators $\overline x_1, \ldots , \overline x_k$ regarding to
circle composition $(\circ)$, $1 - u(x_1, \ldots , x_k) =
\overline u(\overline x_1, \ldots , \overline x_k)$. As $1 -
u(x_1, \ldots ,x_k) \in I(M)$ then by Lemma 5 $\overline
u(\overline x_1, \ldots , \overline x_k) \in I(\overline M)$ and
$\overline u(\overline x_1, \ldots , \overline x_k) = \overline u$
has a normal form. Hence $l(\overline u) > 0$ and, consequently,
$\frak c(\overline u) = 1$. But by Lemma 1 $\frak c(\overline u) =
0$ as $\overline u \in I(\overline M)$. We get a contradiction
with $\frak c(\overline u) = 1$. Hence our supposition that $H
\neq \{1\}$ is false. This completes the proof of Theorem 1.
\smallskip\\

\textbf{Remark.} The proof of Lemma 3 has a constructive character
for free Moufang loops. But Lemma 3 holds for algebras of
$\Omega$-words (see, for example, [17]). Any relatively free
Moufang loop is an algebra of $\Omega$-words. From here it follows
that the Lemma 3 is true for any relatively free Moufang loop.
Then it is easy to see that the main result of this paper (Theorem
1) holds for every relatively free Moufang loop.
\smallskip\\

Further we identify the loop $(\overline M,\star)$ with
$(M,\cdot)$. Then every element in $FM/I(M)$ has the form $\sum_{q
\in M}\lambda_qq$, $\lambda_q \in F$. Further for the alternative
algebra $FM/I(M)$ we use the notation $FM$ and we call them
\textit{''loop algebra''} (in inverted commas). Let $H$ be a
normal subloop of $M$. We denote the ideal of ''loop algebra''
$FM$, generated by the elements $1 - h$ ($h \in H$) by $\omega H$.
If $H = M$, then $\omega M$ will be called the
\textit{''augmentation ideal''} (in inverted commas) of ''loop
algebra'' $FM$. Let us determine the homomorphism $\varphi$ of
$F$-algebra $FM$ by the rule $\varphi(\sum \lambda_qq) = \sum
\lambda_qHq$. By analogy to Lemma 1 it is proved.
\smallskip\\

\textbf{Proposition 1.} \textit{Let $H$ be a normal subloops of
free Moufang loop $M$ and let $FM$ and $\omega M$ are respectively
''loop algebra'' and ''augmentation ideal'' of $M$. Then}

\textit{1) $\omega H \subseteq \text{Ker}\varphi$;}

\textit{2) $1 - h \in \text{Ker}\varphi$ if and only if $h \in
H$;}

\textit{3) $\omega M = \{\sum_{q \in M}\lambda_qq \vert \sum_{q
\in M}\lambda_q
 = 0\}$;}

\textit{4) the ''augmentation ideal'' $\omega M$ is generated as
$F$-module by the elements of the form $1 - q$ ($q \in M$).}
\smallskip\\

Let $\omega M$ denote the augmentation ideal of loop algebra
(without commas) $FM$ and let $\overline{\omega M}$ denote the
''augmentation ideal'' of ''loop algebra''  $FM$.  Then from 4) of
Lemma 1 and 3) of Proposition 1 it follows that
$$\overline{\omega M} = \omega M/I(M).$$

Any Moufang loop $L$ has a representation $L = L/H$, where $L$ is
a free Moufang loop. As we have remarked above, in [15] there are
constructed Moufang loops $L$ that are not embedded into a loop of
invertible elements of any alternative algebras. Then for such
normal subloop $H$ of $L$ $\text{Ker}\varphi = FL$ and by 2) of
Proposition 1 the inclusion $\omega H \subset \text{Ker}\varphi$
is strict.

We mention that Proposition 1 holds also for Moufang loops for
which the Theorem 1 is true.

\smallskip
Tiraspol State University of Moldova$$
$$

The author's home address:

Sandu Nicolae Ion

Deleanu str 1,

Apartment 60

Kishinev MD-2071, Moldova

E-mail: sandumn@yahoo.com

\end{document}